\documentclass[a4paper,twoside,11pt,reqno]{amsart}
\usepackage{hyperref}
\hypersetup{colorlinks=true,citecolor=blue}
\usepackage[latin1]{inputenc} 
\usepackage{color, euscript,graphicx,tabularx,epsfig,psboxit}
\usepackage{graphicx}
\usepackage{eurosym}
\usepackage[latin1]{inputenc}
\usepackage{color}
\usepackage{t1enc}
\usepackage{url}
\usepackage{amsfonts}%
\usepackage{amssymb}%
\usepackage{color,euscript,graphicx,tabularx,epsfig,psboxit,pstricks,amsmath}
\usepackage{graphicx}
\def\h1{\hspace{1cm}}
\def\h2{\hspace{2cm}}\def\h3{\hspace{3cm}}\def\h4{\hspace{4cm}}\def\h5{\hspace{5cm}}
 \definecolor{grey}{rgb}{0.75,0.75,0.75}
\definecolor{orange}{rgb}{1.0,0.5,0.5}
\definecolor{brown}{rgb}{0.5,0.25,0.0}
\definecolor{pink}{rgb}{1.0,0.5,0.5}

\def\dis{\displaystyle}
 
\def\paragraph#1{\textit{#1}.}
\def\C{{\mathbb C}}

\newtheorem{tm}{Theorem}

\newtheorem{lm}{Lemma}

\newtheorem{defi}{Definition}

 \def\cqfd{\unskip\kern 6pt\penalty 500
\raise -2pt\hbox{\vrule\vbox to5pt{\hrule width 4pt
\vfill\hrule}\vrule}}

\def\timederivative#1#2{#1^{(#2)}}

\def\I{\mathbf{I}}
\def\D{\mathbf{D}}

\def\C{\mathbf{C}}
\def\P{\mathbf{P}}
\def\R{\mathbb{R}}

\def\L{\mathbf{L}}

\def\M{\mathbf{M}}

\def\monomials#1{\mathcal{M}}

\newcounter{equationR}
\def\dis{\displaystyle}
\def\Proof{{\bf\noindent Proof. }}

\renewcommand{\thefootnote}

\title[On an Index Reduction Method by Deflation for DAEs]{On an Index Reduction Method by Deflation for Differential-Algebraic Equations $ ^{\ast} $}
\author{F.~Monfreda}
\address{Fabien Monfreda \\ Institut de Math\'ematiques de Toulouse \\ \'equipe MIP \\Bureau 302 \\ Universit\'e Paul Sabatier \\ 118 route de Narbonne \\ 31062 Toulouse Cedex 9\\ France.}
\email{fabien.monfreda@math.univ-toulouse.fr} 
 \author{J.-C.~Yakoubsohn} 
\address{Jean-Claude Yakoubsohn \\ Institut de Math\'ematiques de Toulouse \\ \'equipe MIP\\
Bureau 120 \\ Universit\'e Paul Sabatier \\ 118 route de Narbonne \\ 31062 Toulouse Cedx 9\\
France.}
\email{yak@mip.ups-tlse.fr} 
 
\date{Version of \today}

\begin{document}

\maketitle 
\footnote[1]{\it $ ^{\ast}$This work has been supported by the French ANR-10-BLAN 0109.} 
   
\begin{abstract}
This paper studies a deflation method to reduce and to solve linear differential-algebraic equations (DAEs). It consists to define a sequence of DAEs with index reduction of one unit by step. This is simultaneously performed by substitution and differentiation. At the end of process, we obtain at most an ODE and a list of algebraic constraints which solve the initial DAE. We show on classical examples how works the method. Moreover, we explain how this method extends in the case of linear time-varying DAEs. 
\\\\ {\it AMS classification : 34A09, 65L80.}
\\\\ {\it Keywords :} Differential-algebraic equations ; index reduction ; Kronecker index ; differentiation index.
\end{abstract}

\section{Introduction}
This paper focuses on the study of a deflation type method to solve linear differential-algebraic equations (DAEs). Deflation type methods arise in many problems involving eigenvalues, roots of polynomial systems etc... Roughly speaking, a deflation process replaces the original problem by a problem of smaller size. In the case of DAEs, it consists to determine a sequence of DAEs of size strictly decreasing ; at the end of this process, we obtain at most an ODE and a list of algebraic constraints satisfied by the general solution of the DAE. In this paper, we deal with both linear time-invariant DAEs and linear time-varying DAEs. \\
The literature on linear DAEs is extremely rich with a large classical documentation. In the recent past, the monographs of Kunkel and Mehrmann~\cite{KM06} and Riaza~\cite{riaza08} give an excellent overview and permit to understand the different approaches from the works of pioneers, Weierstrass and Kronecker. The linear case is related to properties of matrix pencils, which are well described in ~\cite{GA259}, where the notion of regularity provides  fundamental results. The diversity of perspectives and techniques for the DAEs is probably due to the central notion of index. In fact, if the notion of Kronecker index in the linear time-invariant DAEs' case is an unifying concept, it is not the same for linear time-dependent case. Several index notions appear according to the point of view adopted : differentiation index developed by Campbell, Gear and Petzold ~\cite{BCP96} ; reduction methods and geometric index studied by Rheinboldt and Rabier~\cite{rh84} and ~\cite{HW10} ; projector-based methods and tractability index introduced by Griepentrog and M\"arz ~\cite{GM86} ; strangeness index by Kunkel and Mehrmann ~\cite{KM06} ; perturbation index by Hairer ~\cite{HW10} ; structural index by Jacobi~\cite{FO09} and Pryce~\cite{Pryce01}. Except for the structural and strangeness indices, all these indices are equal in the linear time-invariant DAEs' case. \\
The complexity of the deflation method, which is proposed here, deals with the Kronecker index (time-invariant case) and the differentiation index (time-varying case). The technical background of the method uses elementary algebra as LU decomposition or singular value decomposition performed on the matrix $ E $ of the following DAE : $ E\dot x=Ax+f. $ \\
In order to provide some motivations at our work, we briefly remember the main contributions in the field of linear DAEs with constant coefficients. \\
The first remark concerns the regularity assumption for linear DAE. In fact, if we don't use this assumption, there are more than countably many different solutions as it is well described in ~\cite{GM86} pages 14-15. The classical reduction of this DAEs type deals with the Weierstrass canonical form and Kronecker index .
Then the linear DAE  is decoupled in two subproblems of type
$\dot u=Ju+g$ and $N\dot v=v+h$, 
where the matrix $J$ and the nilpotent matrix $N$ are related to the Jordan form of $( \lambda E+A)^{-1} E $, for some $ \lambda \in \mathbb{R} $. From a computer algebra point of view, this way permits to obtain
an expression of the solution, i.e $v(t)=-\sum_{k=0}^{\nu-1} N^kh^{(k)}(t)$ where $\nu$
is the Kronecker index, i.e $N^{\nu-1}\ne 0$ and $N^\nu=0$.
But from a numerical analysis point of view, finding the Jordan normal form is known to be numerically unstable : small errors can make large differences in the result, see ~\cite{GW76} or the discussion in ~\cite{GVL96} sec. 7.6.5.  
We also note the Weierstrass canonical form needs to find $\lambda$ such that $\lambda E+A$ is invertible and to compute $(\lambda E+A)^{-1}E$ before to perform a Jordan normal form. These numerical drawbacks led us to propose a reduction of type deflation. \\
On the other hand, the differentiation method, described in  ~\cite{BCP96} page 20, is based on the differentiation of constraints. A new DAE is obtained and the procedure is repeated until to obtain an ODE. The result (Theorem 2.3.3~\cite{BCP96}) says that the number of steps is precisely the Kronecker index in the time-invariant case, and each differentiation of constraints reduces the index by one. But the final ODE has many additional solutions. 
Instead the deflation method proceeds by substitution, and the differentiation is performed during this stage.
 We will state that the number of steps is the minimum of the Kronecker index and the rank. 
 Moreover, solutions of the reduced system are those of the initial system.
 There are other types of reduction in the literature ; mainly the one introduced in
 ~\cite{GM86}, based on projector methods and generalized in \cite{GM89} for DAEs with higher index. A good introduction of these projector methods is done in ~\cite{riaza08} : the solution of the DAE is expressed thanks to matrix chains, which length is equal to the Kronecker index. However, our method appears to be technically simpler. \\
More recently in~\cite{TI08}, an index reduction based on substitution method has been developed for DAEs of the type $E\dot x=Ax+ f(t)$ with
\begin{center}
$ E=\left (\begin{array}{cc}0& K\\L&M\end{array}\right ),\quad A=
 \left (\begin{array}{cc}-B& 0\\0&0\end{array}\right ),\quad f(t)=
 \left (\begin{array}{cc}f_1(t)\\f_2(t)\end{array}\right )
 $\end{center}
where $B$ is an invertible matrix. Introducing the Schur complement $ D = M-LB^{-1}K $,
we obtain a new DAE 
 \begin{equation}\label{DAEGR1}D\dot x_2=f_2-LB^{-1}\dot f_1\end{equation}
 under the constraints
\begin{equation}\label{DAEGR2}Bx_1=f_1-K\dot x_2.\end{equation}
It is proved the index of  the DAE~\ref{DAEGR1} is one lower than that of initial DAE. Next the DAE ~\ref{DAEGR1} is numerically solved and the solution of the initial DAE is recovered thanks to the system of linear equations \ref{DAEGR2}. This method is mainly proposed in the context of electrical circuit where it is possible to find a non-singular constant submatrix $B$ of $\lambda E+A$. Certainly a reduction of an index unit improves the numerical results, but the drawback of this method is the DAE with higher index. In some sense, our deflation method generalizes this approach. Finally, we have to mention the paper of Linh and Mehrmann ~\cite{LM09} where a similar transformation of the DAE and the Schur complement are used in the strangeness-free context. \\ \\
In this paper, we will separately treat the time-invariant case and the time-varying case because of their structural differences.
In section 2, we present the deflation method in the time-invariant case and the main results, which will be proved in section 3. The section 4 gives some examples which illustrate the method. Finally the section 5 shows how works the deflation method in the time-varying case.
\section{Deflation method for linear time-invariant DAEs.}
A linear time-invariant DAE is a problem of the form :
\begin{equation} \label{DAEL} E \dot{x}(t) = A x(t) + f(t), \end{equation}
where $E$ and $A$ are constant matrices in $ \mathbb{R}^{n\times n}$ and $f$ : $I\rightarrow \R^n$ is a sufficiently smooth function, defined on an open interval $I \subseteq \mathbb{R}$. We assume the rank of $E$ is $r$ and this DAE is regular : there exists $\lambda \in \mathbb{R}$ such that the matrix $\lambda E+A$ is non-singular~\cite{GA259}.
\\
Using both basic linear algebra and substitution of certain variables, the main idea is to  separate the differential part and the algebraic part of the DAE ~\ref{DAEL}.
This goal  can be attained using  a decomposition of the matrix as LU decomposition or singular value decomposition (SVD) which are well studied in ~\cite{GVL96}.
For example if  the LU decomposition is used, $E=LU$ and   the DAE ~\ref{DAEL} is equivalent to 
$$U\dot x(t)=L^{-1}Ax(t)+L^{-1}f(t).$$
Hence the algebraic constraints appear since  there are $n-r$ zero rows in the matrix $U$. A similar result is obtained using SVD.  For this reason, we present the deflation
method without reference to the decomposition used to factorize the matrix
$E$. The following lemma is the key of the deflation method.
\begin{lm} \label{regular}
There exist an invertible matrix $U$ and a permutation matrix $P$ such that $\lambda E+ A$ reads
$$
\lambda E+A=U\left (\begin{array}{cc}\lambda S+K&\lambda T+L\\M&N
\end{array}\right )P^{-1}
$$
where $S,K\in\R^{r\times r}$; $T,L\in\R^{r\times n-r}$; $M\in\R^{n-r\times r}$ and $N\in\R^{n-r\times n-r}$ is a non-singular matrix.
\end{lm}   
\Proof
Since the rank of $E$ is $r$, we can factorize $E=U\Sigma:=U\left (\begin{array}{c}F\\0 \end{array}\right )$ where $U$ is invertible and $F$ of size $r\times n$ with rank $r$.
 Writing $U^{-1}A=\left (\begin{array}{cc} A_1\\A_2
\end{array}\right )$ where $A_1\in \R^{r\times n}$ and $A_2\in \R^{n-r\times n}$, we obtain
$$\lambda E+A=U \left (\begin{array}{cc}\lambda F +A_1 \\ A_2 
\end{array}\right ).$$
The matrix $A_2$  has full rank since $\lambda E +A$ is invertible. Hence there exists an invertible submatrix $N\in\R^{n-r\times n-r} $ of $A_2$. Introducing an appropriate permutation matrix $P$, we obtain the desired formula.\cqfd \\ \\
With the notations above, it is easy to see the DAE \ref{DAEL} is equivalent to
\begin{equation}\label{DAEL2}
 \left (\begin{array}{cc} S&T \\ 0&0 
\end{array}\right )P^{-1} \dot x=  \left (\begin{array}{cc} K&L \\ M&N 
\end{array}\right )P^{-1} x+ \left (\begin{array}{cc} g \\h 
\end{array}\right ),
\end{equation}  
where $U^{-1}f=\left (\begin{array}{c}g\\h\end{array}\right )$. 
After an easy computation, the DAE \ref{DAEL} is finally reduced to
\begin{equation}\label{DAEL5}
(S-TN^{-1}M)\dot u=(K-LN^{-1}M)u +TN^{-1}\dot h-LN^{-1}h +g 
\end{equation}
\begin{equation}\label{DAEL6}
v=-N^{-1}M u -N^{-1}h
\end{equation}
\begin{equation}\label{DAEL7}
P^{-1}x=(u,v)^T
\end{equation} 
The size of the DAE~\ref{DAEL5} is less than that of the DAE ~\ref{DAEL}.
This conduces to define a deflated DAE of an initial DAE.
\begin{defi}\label{deflation}
We note \begin{align*}
&E_1=S-TN^{-1}M,\\ 
&A_1=K-LN^{-1}M, \\
&f_1=TN^{-1}\dot h-LN^{-1}h +g,\\&x^1=u.
\end{align*}
We say that  $E_1\dot x^1=A_1x^1+f_1$ is a deflated DAE of \ref{DAEL}.
Moreover $P$ is the permutation matrix  and
$Mu+Nv+h=0$ is the algebraic constraint associated to this deflation.
\end{defi}
Let us remark there are in general several deflated DAEs of a given DAE~: in fact  it depends upon the choice of matrices $M$ and $N$. Fortunately the ranks of matrices $E_1$ and $A_1$ are invariant as it is stated in the following.
\begin{tm}\label{equivalence}
Let us consider two deflated DAEs of the DAE ~\ref{DAEL} represented respectively
by the matrix pencils $(E_1,A_1)$ and $(\tilde E_1,\tilde A_1)$. Then we have $rank(E_1)=rank(\tilde E_1)$
and $rank(A_1)=rank(\tilde A_1)$.
 \end{tm}
 Moreover, the regularity is preserved during a step of  deflation.
 \begin{tm}\label{regularity}
 If $\lambda E+A$ is non-singular then $\lambda E_1+A_1$ is non-singular.
 \end{tm}   
 From the two previous theorems which will be proved in the next section, it follows first that the deflation process is independent
 of the choice of the two matrices $M$ and $N$, and next that the regularity is preserved. This suggests the following reduction procedure.
\\ \\
{\bf Deflation algorithm.}
\\ \\
{\bf Input }:   $E_0=E$, $A_0=A$, $f_0=f$, $x^0=x$, $r_{-1}=n$.
\\ {\bf Step $j+1$, $j\ge 0$.}\\
   If  $E_{j}$ is singular and $E_{j}\ne 0$ 
   \begin{itemize}
\item[1--] Let  $r_{j}$ the rank of $E_j$.
\item[2--] Compute  $(E_{j+1},A_{j+1},f_{j+1})$ from $(E_{j},A_{j},f_{j})$ using the formulas of the definition ~\ref{deflation}. Let $P_j$ the permutation matrix of this reduction. 
\item[3--] Compute the change of variable $(x_1^{j},x_2^{j})^T:=P_{j}^{-1}x^j$ where the size of $x_1^{j}$ is $r_{j}$.
\item[4--]  Compute the algebraic constraint : $0=M_jx_1^{j}+N_jx_2^{j}+f_{j,r_{j}+1:r_{j-1}}$, where $ f_{j,r_{j}+1:r_{j-1}} $ means the  coordinates  $r_j+1$ to $r_{j-1} $ of the vector $ f_{j} $. 
\item[5--] Let $x^{j+1}:=x_1^{j}$.
\end{itemize}
else stop.\\
{\bf Output} For $j\ge 0$ the sequences 
$$\begin{array}{cc}
\textrm{DAEs}&\textrm{algebraic constraints}\\
E_{j+1}\dot x^{j+1}=A_{j+1}x^{j+1}+f_{j+1},&0=M_jx_1^{j}+N_jx_2^{j}+f_{j,r_{j}+1:r_{j-1}}.
 \end{array}$$
The main goal of this paper is to prove the deflation algorithm stops in a finite number of steps. More precisely 
\begin{tm}\label{reduction_algo}
The number of steps of the deflation algorithm is bounded by $min(rank(E),ind(E,A))$
where $ind(E,A)$ is the Kronecker index of the pencil $(E,A)$.
Moreover, the ranks of matrices $E_j$ and $A_j$ determine a sequence of invariants which are characteristic for $E$ and $A$. 
\end{tm}
The general solution is described by the following result.
\begin{tm}\label{solution}
Let $k$ be the number of steps of the deflation algorithm. Then the coordinates of the solution satisfy
\begin{equation}\label{sol_DAE1}
E_k\dot x^{k}=A_kx^k+f_k
\end{equation}
where $E_k$ is invertible or equal to zero and 
\begin{equation}\label{sol_DAE2}
x^j_2=-N_j^{-1}M_jx^j_1-N_j^{-1}f_{j,r_{j}+1:r_{j-1}},\quad 0\le j\le k-1.
\end{equation}
\end{tm}
\Proof From the definition of the deflation algorithm.
\section{Proofs of theorems of the section 1. }
To prove that let us remember some fundamental notions.
The index of a matrix $B$ is the smallest integer which verifies $Ker(B^k)=Ker(B^{k+1})$.
If the index $B$ is equal to $k$ then it is equivalent to the transversality condition or the range-nullspace decomposition ~\cite{Meyer00}
page  394:
$$\R^n=Ker(B^k)\oplus Im(B^k).$$
If the index of $B$ is zero then $B$ is non-singular.
The Kronecker index of the DAE \ref{DAEL} is the index of the matrix $(\lambda E+A)^{-1}E$.
We denote by  $ind(E,A)$ this index. Moreover, lemma 7 page 196 of ~\cite{GM86} shows that
$ ind(E, A) = ind(EQ, AQ) $, if $ Q $ is invertible. \\
We first prove the theorem ~\ref{equivalence}. \\ \\
{\bf Proof of theorem~\ref{equivalence}.}
Let us suppose there exist two choices of matrices
$(S,T,K,L,M,N)$ and $(\tilde S,\tilde T,\tilde K,\tilde L,\tilde M,\tilde N)$
for which the lemma~\ref{regular} holds.
Then there exists a permutation matrix $P_0$ such that
$$
\left (\begin{array}{cc}S&T\\M&N\end{array}\right)P_0=
\left (\begin{array}{cc}\tilde S&\tilde T\\\tilde M&\tilde N\end{array}\right ).
$$
Since $N$ and $\tilde N$ are non-singular, we can write
\begin{align*}
&\left (\begin{array}{cc}TN^{-1}&I\\I&0\end{array}\right )
\left (\begin{array}{cc}N&0\\0&E_1\end{array}\right )
\left (\begin{array}{cc}N^{-1}M&I\\I&0\end{array}\right )
P_0\\
&\quad\quad\quad\quad\quad\quad\quad
=
\left (\begin{array}{cc}\tilde
T\tilde N^{-1}&I\\I&0\end{array}\right )
\left (\begin{array}{cc}\tilde N&0\\0&\tilde E_1\end{array}\right )
\left (\begin{array}{cc}\tilde N^{-1}\tilde M&I\\I&0\end{array}\right )
\end{align*}
Hence the ranks of the matrices $E_1$ and $\tilde E_1$ are equal.
By a same way, we show $rank(A_1)=rank(\tilde A_1)$.\cqfd
\\
 We now prove the theorem ~\ref{regularity}.
\\ \\
{\bf Proof of theorem~\ref{regularity}.}
From lemma 1 it follows the matrix $\lambda \Sigma P+U^{-1}AP$ is non-singular.
Using the Schur complement of this matrix, we have
\begin{align*}
\lambda \Sigma P+U^{-1}AP&= \left (\begin{array}{cc}
\lambda S+K&\lambda T + L\\ M&N 
\end{array}\right )\\
&=\left (\begin{array}{cc}
(\lambda T+L)N^{-1}&I\\I&0
\end{array}\right )
\left (\begin{array}{cc}
N&0\\0&\lambda E_1+A_1
\end{array}\right )
\left (\begin{array}{cc}
N^{-1}M&I\\I&0
\end{array}\right )
\end{align*}
It follows the matrix $\lambda E_1+A_1$ is non-singular. \cqfd
\\
To prove the theorem ~\ref{reduction_algo} we need some lemmas.
\begin{lm}\label{inverse}
With the notations of the section 1 and $E=U\Sigma$, $C_1=\lambda E_1+A_1$, we have :
\quad\begin{itemize}
\item[1--] $(\lambda \Sigma P+U^{-1}A P)^{-1}\Sigma P=\left (\begin{array}{cc}
C_1^{-1}S&C_1^{-1}T\\-N^{-1}MC_1^{-1}S&-N^{-1}MC_1^{-1}T
\end{array}\right )$
\item[2--] 
\begin{align*}
\left ((\lambda \Sigma P\right .&\left . +U^{-1}AP)^{-1}\Sigma P\right )^k\\
&=\left (\begin{array}{cc}
(C_1^{-1}E_1)^{k-1}C_1^{-1}S&(C_1^{-1}E_1)^{k-1}C_1^{-1}T\\
-N^{-1}M(C_1^{-1}E_1)^{k-1}C_1^{-1}S&-N^{-1}M(C_1^{-1}E_1)^{k-1}C_1^{-1}T
\end{array}\right )
\end{align*}
\end{itemize}
\end{lm}
{\Proof} A straightforward computation gives the result of the part 1. \\
For the part 2, we apply the lemma below to the identity of the part 1.\cqfd
\begin{lm} Let $k\ge 1$.
$$\left (\begin{array}{cc}A&B\\CA&CB\end{array}\right )^k=
\left (\begin{array}{cc}(A+BC)^{k-1}A&(A+BC)^{k-1}B\\
C(A+BC)^{k-1}A&C(A+BC)^{k-1}B\end{array}\right )$$
\end{lm}
{\Proof} By induction.\cqfd \\ \\
We next state a result which appears in ~\cite{GM86}, in theorem 13 page 198. We remember the proof for sake of completion.

\begin{lm}\label{index1}The index of the DAE \ref{DAEL} is equal to one iff
 $x\in Ker\,E$ and $Ax\in Im\, E\Rightarrow x=0.$ 
 \end{lm}
{\Proof} Let us suppose the index of the DAE \ref{DAEL} is one. This implies 
$\R^n=Im(\lambda E+A)^{-1}E\oplus Ker(\lambda E+A)^{-1}E.$
Let us suppose $x\in Ker\, E$ and $Ax \in Im\, E.$ Then $x\in Ker\, E$ implies $x\in Ker\, (\lambda E+A)^{-1}E$. 
On the other hand, there exists $y\in\R^n$
 such that $Ax=Ey$. Since $(\lambda E+A)x=Ax$ and $\lambda E+A$ is non-singular we have
 $x=(\lambda E+A)^{-1}Ey\in Im\,(\lambda E+A)^{-1}E$. 
Hence $x\in Ker (\lambda E+A)^{-1}E\cap Im (\lambda E+A)^{-1}E$. 
It follows $x=0$.
 \\
 Conversely, let us suppose $x\in Ker\,E$ and $Ax\in Im\, E\Rightarrow x=0$. Assuming the index $k>1$, there exists $y\ne 0$ such that $(\lambda E+A)^{-1}Ey\ne 0$
 and $\left( (\lambda E+A)^{-1}E \right)^2y=0$.  Hence $z:=(\lambda E+A)^{-1}Ey\in  Ker\, E$. Prove that $Az\in Im\,E$. In fact $(\lambda E+ A)z=Ey$, hence $Az=Ey\in Im\,E$. This means $z=0$
 in contradiction with the definition of $z$. 
  \cqfd
  \begin{lm}\label{nonsingular}The index of the DAE \ref{DAEL} is equal to one iff
  $E_1$ is non-singular.
 \end{lm}
{\Proof} Let us suppose $ind(E,A)=1$ and prove that $Ker E_1=\{0\}$. It is equivalent to show that the
matrix $\dis \left (\begin{array}{cc} S&T\\M&N\end{array}\right )$ is non-singular since
 the Schur complement of this matrix is $E_1$.
Let us consider $x$ such that
$\dis \left (\begin{array}{cc}S&T\\M&N\end{array}\right )x=0$.  
It implies both $x\in Ker(S,T)=Ker (\Sigma P)$ and $(M,N)x=0$.  Hence
$\dis \left (\begin{array}{cc}K&L\\M&N\end{array}\right )x=
 \left (\begin{array}{c}(K,L)x\\0\end{array}\right )$.
 Since the rank of $\Sigma P$ is $r$, there exists $y\in\R^n$ such that
 $\left (\begin{array}{c}(K,L)x\\0\end{array}\right )=\Sigma P y$.
 In fact we have $x\in Ker (\Sigma P)$ and $U^{-1}APx=\Sigma P y\in Im(\Sigma P)$. From lemma \ref{index1} it follows $x=0$ and the matrix $\dis \left (\begin{array}{cc} S&T\\M&N\end{array}\right )$ is non-singular.
 \\
 Let us suppose now $E_1$ is non-singular.\\ From lemma \ref{regular}, we have
 $ind(E,A)=ind(\Sigma P, U^{-1}AP)$. From lemma \ref{index1}, proving that $ind(\Sigma P, U^{-1}AP)=1$ is equivalent to establish the assertion $x\in Ker(\Sigma P)$ and 
 $U^{-1}APx\in Im(\Sigma P) \Rightarrow x=0$. If we have $x\in Ker(\Sigma P)$ and 
 $U^{-1}APx\in Im(\Sigma P)$, this implies $(S,T)x=0$ and $(M,N)x=0$. Since 
 $\dis \left (\begin{array}{cc} S&T\\M&N\end{array}\right )$ is non-singular, it follows $x=0$.
 \cqfd
\\
We can now state the result which links the indices of the DAE~\ref{DAEL} and  the deflated DAE~\ref{DAEL5}.
\begin{tm}\label{index_reduction}
If $E_1\ne 0$ then $ind(E_1,A_1)=ind(E,A)-1.$ 
\end{tm} 
{\Proof}
The case $ind(E,A)=1$ is treated by the lemma \ref{nonsingular}.\\
Let us suppose first that $k=ind(E,A)\ge 2$.
Since $(\lambda E+A)^{-1}E=P\,(\lambda \Sigma P + U^{-1}AP)^{-1}\Sigma PP^{-1}$,
the indices $ind(E,A)$ and $ind(\Sigma P, U^{-1}AP)$ are equal. We then have
\begin{equation}\label{directsum}
\R^n=Ker\left ((\lambda \Sigma P + U^{-1}AP)^{-1}\Sigma P\right )^k\oplus Im \left( (\lambda \Sigma P + U^{-1}AP)^{-1}\Sigma P \right)^k.
\end{equation}
It is sufficient to show 
$$x\in Ker\, (C_1^{-1}E_1)^{k-1}\cap
 Im\,(C_1^{-1}E_1)^{k-1}\Rightarrow x=0,$$
where $C_1=\lambda E_1+A_1$.
 If $x\in Ker\,(C_1^{-1}E_1)^{k-1}$ then $(C_1^{-1}E_1)^{k-2}C_1^{-1}(Sx-TN^{-1}Mx)=0$. 
 From the identity of lemma \ref{inverse} part 2, it follows $$(x,-N^{-1}Mx)^T\in 
Ker \left( (\lambda \Sigma P + U^{-1}AP)^{-1}\Sigma P\right )^k.$$
 Now if  $x\in Im\,(C_1^{-1}E_1)^{k-1}$, there exits $y\in\R^r$
such that
 $x=(C_1^{-1}E_1)^{k-1}y$. 
Since the rank of $\Sigma P$ is equal to $r$,
 there exists $(u,v)^T\in\R^r\times\R^{n-r}$ 
 such that $y=C_1^{-1}(Su+Tv)$. Then  $x=(C_1^{-1}E_1)^{k-1}C_1^{-1}(Su+Tv)$. 
Always from lemma \ref{inverse} part 2, it follows 
$$(x,-N^{-1}Mx)^T\in 
Im\, \left ((\lambda \Sigma P + U^{-1}AP)^{-1}\Sigma P\right )^k.$$
From \ref{directsum} we deduce $(x,-N^{-1}Mx)^T=0$ and finally $x=0$.\cqfd
 \\
Next, we state some properties concerning the rank of $E_1$.
\begin{lm}\label{rankE}
Let us consider the regular DAE \ref{DAEL}.
\begin{itemize}
\item[1--]  If  $ind(E,A)>1$ then
$rank(E)>rank(E_1)$.
 \item[2--]  If  $ind(E,A)=1$ then
$rank(E)=rank(E_1)$.
\item[3--] If $ind(E,A)>1$ and $rank(E)=1$ then $E_1=0$.
\end{itemize}
\end{lm}
{\Proof} Prove the part 1. From construction $r=rank(E)\ge rank(E_1)$ since $E_1\in\R^{r\times r}$. If $rank(E_1)=r$, this implies $E_1$ is non-singular  and we have $ind(E,A)=1$
 from lemma  \ref{nonsingular}. This contradicts the hypothesis $ind(E,A)>1$.
 The part 2 is a direct consequence of part 1. For the part 3, we observe that $E_1\in \R$ since $rank(E)=1$. If $E_1\ne 0$ then $E_1^{-1}$ exists. This contradicts $ind(E,A)>1$. Hence $E_1=0$.\cqfd
\\ Now, we are able to prove the theorem \ref{reduction_algo}. \\ \\
{\bf Proof of theorem \ref{reduction_algo}.}
At each step of the algorithm, the index is strictly decreasing from theorem \ref{index_reduction}. In the same way, from lemma \ref{rankE}, the rank of each $E_j$ is also strictly decreasing. The reduction algorithm stops when the index or the rank is $0$. The result follows.\cqfd 
\section{Examples}
We present in the  examples below, the different steps of the deflation algorithm with the differential part in the left
and, the algebraic constraints in the right. \\ \\
{\bf Example 1.}~\cite{BCP96} page 19.
$$\left (\begin{array}{ccc}1&0&0\\0&1&0\\0&0&0\end{array}\right )\dot x=
\left (\begin{array}{ccc}0&0&-1\\-1&0&0\\0&-1&0\end{array}\right )x +
 f(t).$$
The Kronecker index is $2$.   
 The sequence of DAEs and the algebraic constraints given by the deflation algorithm is successively 
 described in the two steps below. \\ \\
 { \it Step 1.}
 The permutation matrix consist to swap the columns 2 and 3. Hence $S-TN^{-1}M=\left (\begin{array}{cc}1&0\\0&0\end{array}\right )-
 \left (\begin{array}{c}0\\1\end{array}\right )(-1)\left (\begin{array}{cc}0&0\end{array}\right )=S.$
 The deflated DAE and the constraints are :
 $$\left (\begin{array}{cc}1&0\\0&0\end{array}\right )
 \left (\begin{array}{c}\dot x_1\\\dot x_3\end{array}\right )=
\left (\begin{array}{cc}0&-1\\-1&0\end{array}\right )x+
  \left (\begin{array}{c}f_1\\f_2-\dot f_3\end{array}\right ),\quad\quad 0=x_2-f_3$$
{\it Step 2.}
  $$0=-x_3+f_1-\dot f_2+\ddot f_3,\quad\quad 0=x_1-f_2+\dot f_3.$$
  In this case the algorithm directly gives algebraic constraints which determine the solution.
  \\ \\
{\bf Example 2.} $N_k\dot x=x$ where $N$ is an elementary nilpotent matrix of size $ k \times k $ :  
$N_k=\left (\begin{array}{ccccc}0&1\\&\ddots&\ddots\\
&&0&1\\&&&0\end{array}\right )$.
We have $N^{k-1}\ne 0$ and $N^k=0$. 
      \\The deflated  DAE   is $N_{k-1}\dot x_{1:k-1}=x_{1:k-1}$ and the 
constraint relation $x_k=0$.
\\ \\
{\bf Example 3.} This class of DAEs appears in~\cite{RMB00} and has a Kronecker index equal to $1$.
\begin{align*}
\dot x_2+\dot x_3&=-x_1+f_1(t)\\
\dot x_2+\dot x_3&=-x_2+f_2(t)\\
\dot x_4+\dot x_5&=-x_3+f_3(t)\\
\dot x_4+\dot x_5&=-x_4+f_4(t)\\
0&=-x_5+f_5(t)
\end{align*}
The deflation algorithm stops after one step. \\ \\
{\it Step 1.}
$$
\begin{array}{c}
\dot x_2+\dot x_4=-x_2+f_2-\dot f_3+\dot f_4\\
\dot x_4=-x_4+f_4-\dot f_5
\end{array}\quad\quad
0=\left (\begin{array}{c}
x_5-f_5\\
x_3-x_4-f_3+f_4\\
x_1-x_2-f_1+f_2
\end{array}
\right )
$$
{\bf Example 4.} This  is the example 7 ~\cite{TI08} which described an electric circuit with index three provided by ~\cite{GR96}. This DAE is first reduced to a DAE with index two and next
the solution is numerically computed thanks to the DAE solver in Matlab.
For this example, our method finds an ODE which gives the exact solution.
In our context this DAE reads
$$\begin{array}{c}
0=-x_1-x_2\\
0=-x_5+x_6\\
0=-x_4+V(t)\\
C\dot x_8=x_2\\
L\dot x_7=x_6 \\
0=-ax_1-x_3\\
\dot x_3+\dot x_7=0\\
\dot x_4-\dot x_8=0
\end{array}
 $$
The index of this DAE is equal to $3$. The sequence of deflated DAEs given by the algorithm is given in the three steps below. \\ \\
{\it Step 1. }
$$\begin{array}{c}
\dot x_7-a\dot x_1=0\\
-\dot x_8=-\dot V(t)\\
L\dot x_7=x_6\\
C\dot x_8=-x_1\\
L\dot x_7=x_6 \\
 \end{array}
\quad\quad
0=\left (\begin{array}{c}
 x_5-x_6\\
 x_1+x_2\\
 ax_1+x_3\\
 x_4-V(t)
 \end{array}\right )
 $$
{\it Step 2. }
 $$\begin{array}{c}
-\dot x_8=-\dot V(t)\\
\dot x_7=-aC\ddot V(t)\\
0=x_6+LaC\ddot V(t) 
 \end{array}
\quad\quad
0=x_1+C\dot V(t) 
 $$
{\it Step 3. }
 $$\begin{array}{c}
-\dot x_8=-\dot V(t)\\
\dot x_7=-aC\ddot V(t) 
 \end{array}
\quad\quad
0=x_6+LaC\ddot V(t) 
 $$
 Consequently the solution satisfies
 \begin{align*}
 &-\dot x_8=-\dot V(t)\\
&\dot x_7=-aC\ddot V(t) \\
&0=x_6+LaC\ddot V(t) \\
&0=x_1+C\dot V(t) \\
 &0=x_5-x_6\\
 &0=x_1+x_2\\
 &0=ax_1+x_3\\
 &0=x_4-V(t)
 \end{align*}
 \section{ Deflation method for linear time-varying DAEs}
 Now, we consider a linear time-varying DAE
 \begin{equation}\label{DAELTV}
 E(t)\dot x(t) = A(t)x(t) + f(t)
 \end{equation}
 where $E(t)$, $A(t)$, $f(t)$ are matrices (resp. vector) in $\R^{n\times n}$ (resp. $\R^n$) , sufficiently smooth, defined on an open real interval $I$. 
We assume the rank of $E(t)$ is equal to $r$ on the interval $I$. 
The assumption of regularity
of the DAE ~\ref{DAELTV} on the interval $I$ defined in the section 1 does not apply in this case.
The classical example ~\cite{KM06} page 56
$$\left (\begin{array}{cc}
1 &-t\\ 0&0 
 \end{array}\right )\dot x=
\left (\begin{array}{cc}
0&0 \\-1&t
 \end{array}\right ) x+f(t)
 $$
shows there exists a solution of this equation whereas the determinant of $\lambda E+A$ is zero for all $\lambda$.
For this reason, we precise the notion of regularity that we will use. It has been introduced   under another form in definition 3.1 of ~\cite{RARH96}.
\begin{defi}
The DAE ~\ref{DAELTV} is geometrically regular on the open interval $I$ if
\begin{itemize}
\item[1--] the rank of $E(t)$ is constant on $I$.
\item[2--] the conclusion of lemma~\ref{regular} holds. More precisely, there exist a permutation matrix $P$ and matrices $S(t),T(t),K(t),L(t),M(t),N(t)$ sufficiently smooth on $I$
such  that $N(t)$ is invertible and
$$E(t)=U(t)\left (\begin{array}{cc}S(t)&T(t)\\0&0\end{array}
\right )P,\quad
U(t)^{-1}A(t)P=\left (\begin{array}{cc}K(t)&L(t)\\M(t)&N(t)\end{array}\right ).$$
\end{itemize}
\end{defi}
The DAE of the example above is geometrically regular on $\R_+$.
 The following definition specifies the formulas for one step of the deflation process.
 \begin{defi}
Let us suppose the DAE ~\ref{DAELTV} is geometrically regular on the interval $I$.
With the notations of the section 1 and $U(t)^{-1}f(t)=(g(t),h(t))^T$ we define :
\begin{align*}
E_1&=S-TN^{-1}M\\
A_1&=K-LN^{-1}M+TN^{-1}(\dot M-\dot N N^{-1} M)\\
f_1&=TN^{-1}(\dot h-\dot NN^{-1}h)-LN^{-1}h+g \\
P^{-1}x&=(u,v)^T,\quad x^1=u
\end{align*}
where all the matrices and vectors above depend on $t$, unless the matrix $P$.
We say that  $E_1(t)\dot x^1=A_1(t)x^1+f_1(t)$ is a deflated DAE of $E(t)\dot x=A(t)x+f(t)$.
Moreover $P$ is the permutation matrix  and
$Mu+Nv+h=0$ is the algebraic constraint of this deflation.
 \end{defi}
Unlike the case of linear time-invariant DAEs, the deflated  DAE is not necessarily geometrically regular.
For example if 
$$E(t)=\left (\begin{array}{ccc}1&0&1\\1&1&1\\0&0&0\end{array}
\right ),\quad
A(t)=\left (\begin{array}{ccc}0 &0&0\\0&0&0\\t&t&t\end{array}\right ),$$
then
$$E_1(t)=\left (\begin{array}{ccc}0&1\\0&0\end{array}
\right ),\quad
A_1(t)=0.$$
Consequently we need to suppose the regularity of the deflated DAE to continue the deflation
process. In this way the deflation algorithm defined for linear time-invariant DAE is transposable for linear time-varying DAE.
Then, we have
\begin{tm}
Let $k$ be the number of steps  of the deflation algorithm applied to the
linear time-varying DAEs' case. 
Suppose that each DAE $E_j(t)\dot x^j=A_j(t)x^j + f_j(t)$, $0\le j\le k-1$ obtained
in the deflation algorithm is geometrically regular on $I$. Then the differentiation index of the DAE~\ref{DAELTV} is  equal to $k$.
Moreover this index is bounded by the rank of $E(t)$ and the ranks of matrices $E_j(t)$ and $A_j(t)$ determine a sequence of invariants which are characteristic for $E(t)$ and $A(t)$. 
\end{tm}
{\Proof} From the definition of the differentiation index given in ~\cite{BCP96} page 22, it follows this index is equal to $k$. The other properties follow of the definition of the deflation process.\cqfd
\\ \\
 {\bf Example 5.} This example is treated in~\cite{riaza08} page 65. We suppose the functions $C_1(t)$, $C_2(t)$ and $L(t)$ do not vanish.  
\begin{align*}
\frac{d}{dt}(C_1(t)x_1)&=x_4-x_5\\
\frac{d}{dt}(C_2(t) x_2)&=-x_3-x_4\\
\frac{d}{dt}(L(t) x_3)&=x_2\\
0&=x_1-x_2+R_1(t)x_4\\
0&=x_1-R_2(t)x_5
\end{align*}
We note by $E$ and $A$ the matrices of this DAE which the index depends to the  values of the parameters. \\ \\
{\it Index one.} $R_1(t)\ne 0$, $R_2(t)\ne 0$. \\ 
{\it Step 1.}
We permute  the columns $1$ and $4$ and the columns $2$ and $5$ in the matrices $E$ and $A$ simultaneously. 
\begin{align*}
C_1R_2\dot x_5 &=x_4-(1+\dot C_1R_2+C_1\dot R_2)x_5\\
C_2R_1\dot x_4+C_2R_2\dot x_5
&=-x_3-(1+\dot C_2R_1+C_2\dot R_1)x_4-(\dot C_2R_2+C_2\dot R_2) x_5\\
L\dot x_3 &=-\dot Lx_3+R_1x_4 +R_2x_5
\end{align*}
$$
0=\left (\begin{array}{c}
 R_1x_4+x_1-x_2\\
 -R_2x_5+x_1\end{array}\right )
 $$
At this step we have $x^1:=(x_4,x_5,x_3)^T$ and $x^2=(x_1,x_2)^T$. \\ \\
 {\it Index  two.} $R_1(t)= 0$, $R_2(t)\ne 0$. Let us also suppose $C_1+C_2\ne 0$. 
\\ {\it Step 1. }
We permute  the columns $1$ and $4$ and the columns $2$ and $5$ in the matrices $E$ and $A$ simultaneously. 
$$\begin{array}{c}
C_1R_2\dot x_5=x_4-(1+\dot C_1R_2+C_1\dot R_2)x_5\\
C_2R_2\dot x_5=-x_4-(\dot C_2R_2+C_2\dot R_2)x_5-x_3\\
L\dot x_3= R_2x_5-\dot L x_3 
\end{array}\quad\quad
0=\left (\begin{array}{c}
 x_1-x_2\\
 -R_2x_5+x_1\end{array}\right )
 $$
We note $E_1$ and $A_1$ the matrices of this DAE and
 we have $x^1:=(x_4,x_5,x_3)^T$ and $x^2=(x_1,x_2)^T$.
 \\
{\it Step 2. }
We note $E_1:=PE_1$ and $A_1:=PA_1$ with
$P=\left (\begin{array}{ccc}1&-C_1/C_2&0\\0&1&0\\0&0&1\end{array}\right )$.
We permute the rows $1$ and $3$ in the new matrices $E_1$ and $A_1$.
\begin{align*}
L(1+c)c\dot x_4 +Lbc\dot x_5&=-(\dot L (c+c^2)-L\dot c)x_4-(R_2c^2+\dot Lcb+Lc\dot b-L\dot cb)x_5\\
C_1R_2\dot x_5&=x_4+(ca+b)x_5\\
0&= (1+c)x_4+bx_5+cx_3.
\end{align*}
with $c:=C_1/C_2$, $a:=-C_2\dot R_2-\dot C_2R_2$ and  $b:=-1-\dot C_1R_2+\dot C_2R_2c$. 
The matrix $E_2$ is invertible since we have  $C_1(t)+C_2(t)\ne 0$ from the assumption. 
Finally the solution satisfies
\begin{align*}
L(1+c)c\dot x_4 +Lbc\dot x_5&=-(\dot L (c+c^2)-L\dot c)x_4-(R_2c^2+\dot Lcb+Lc\dot b-L\dot cb)x_5\\
C_1R_2\dot x_5&=x_4+(ca+b)x_5\\
0&= (1+c)x_4+bx_5+cx_3\\
0&=x_1-x_2\\
0&= -R_2x_5+x_1.
\end{align*} \\ 
 {\it Index three.} $R_1(t)= 0$, $R_2(t)\ne 0$ and $C_1(t)+C_2(t)=0$.
 Hence $\dot C_1(t)+\dot C_2(t)=0$.
 \\
 {\it Step 1. }
It is the same as step one for index two.
$$\begin{array}{c}
C_1R_2\dot x_5=x_4-(1+\dot C_1R_2+C_1\dot R_2)x_5\\
C_1R_2\dot x_5=x_4-(\dot C_1R_2+C_1\dot R_2)x_5-x_3\\
L\dot x_3= R_2x_5-\dot L x_3 
\end{array}\quad\quad
0=\left (\begin{array}{c}
 x_1-x_2\\
 -R_2x_5+x_1\end{array}\right )
 $$
We note $E_1$ and $A_1$ the matrices of this DAE and
 we have $x^1:=(x_4,x_5,x_3)^T$ and $x^2:=(x_1,x_2)^T$.
\\
 {\it Step 2. } We proceed as in step $2$  for index 2.  
\begin{align*}
L\dot x_5&=-(R_2+\dot L)x_5\\
C_1R_2\dot x_5&=x_4-(1+\dot C_1R_2+C_1\dot R_2)x_5\\
0&= x_5+x_3.
\end{align*}
Now, the matrix $E_2$ is not invertible. 
We have $x^1:=(x_4,x_5)^T$ and $x^2:=x_3$.
\\
{\it Step 3. }
We note $E_2:=PE_2$ and $A_2:=PA_2$ with
$P=\left (\begin{array}{cc}1&-L/(C_1R_2)\\0&1\end{array}\right )$.
We next permute the rows $1$ and $2$ and columns $1$ and $2$ in the new matrices $E_2$ and $A_2$.
Then 
 \begin{align*}
LC_1R_2\dot x_5&=-(La+b)x_5 \\
 0& = bx_5+Lx_4
 \end{align*}
with 
$a=1+\dot C_1R_2+C_1\dot R_2$ and $b=C_1R_2^2+C_1R_2\dot L-L-\dot C_1R_2L-C_1\dot R_2L$.
Finally the solution satisfies
\begin{align*}
LC_1R_2\dot x_5&=-(La+b)x_5 \\
 0& = bx_5+Lx_4\\
 0&= x_5+x_3\\
 0&= x_1-x_2\\
0&= -R_2x_5+x_1.
 \end{align*}
\section{Conclusion}
In this paper, a deflation algorithm for DAEs has been studied. 
It highlights the key role of index notions, especially the Kronecker index and the differentiation index. Moreover, the construction of the algorithm establishes a connection between the rank of $ E $ and the index of the DAE. There is the choice to use LU decomposition or SVD at each deflation step.
The invariants of this algorithm are the index of the DAE and the successive ranks $r_i$ of deflated DAEs. The arithmetic complexity of the deflated sequence is in  
$ \displaystyle O(\sum_{i=-1}^kr_i^3) $ operations.  
This reduction provides at most an ODE, which is solvable with classical techniques, and also algebraic equations. This algorithm has the advantage of being technically simple and therefore it brings an additional process in order to reduce and to solve linear differential-algebraic equations.
\bibliographystyle{plain}
\bibliography{Article_Deflation_Derniere_version}

\begin{thebibliography}{10}

\bibitem{BCP96}
K.E. Brenan, S.L. Campbell, and L.R. Petzold.
\newblock {\em {Numerical solution of initial-value problems in
  differential-algebraic equations}}.
\newblock Society for Industrial Mathematics, 1996.

\bibitem{GA259}
F.R. Gantmacher.
\newblock {\em {The Theory of Matrices: Vol.: 2}}.
\newblock Chelsea publishing company, 1959.

\bibitem{GVL96}
G.H. Golub and C.F. Van~Loan.
\newblock {\em {Matrix computations}}.
\newblock Johns Hopkins Univ Pr, 1996.

\bibitem{GW76}
G.H. Golub and J.H. Wilkinson.
\newblock {Ill-conditioned eigensystems and the computation of the Jordan
  canonical form}.
\newblock {\em SIAM review}, 18(4):578--619, 1976.

\bibitem{GM86}
E.~Griepentrog and R.~M{\"a}rz.
\newblock {\em {Differential-algebraic equations and their numerical
  treatment}}.
\newblock BSB Teubner, 1986.

\bibitem{GM89}
E.~Griepentrog and R.~M{\"a}rz.
\newblock {Basic properties of some differential-algebraic equations}.
\newblock {\em Z. Anal. Anwendungen}, 8(1):25--41, 1989.

\bibitem{GR96}
M.~G{\"u}nther and P.~Rentrop.
\newblock {The differential-algebraic index concept in electric circuit
  simulation}.
\newblock {\em Zeitschrift f{\"u}r angewandte Mathematik und Mechanik},
  76:91--94, 1996.

\bibitem{HW10}
E.~Hairer and G.~Wanner.
\newblock {\em {Solving ordinary differential equations II: Stiff and
  differential-algebraic problems}}.
\newblock Springer, 2010.

\bibitem{KM06}
P.~Kunkel and V.L. Mehrmann.
\newblock {\em {Differential-algebraic equations: analysis and numerical
  solution}}.
\newblock European Mathematical Society, 2006.

\bibitem{LM09}
V.H. Linh and V.~Mehrmann.
\newblock {Lyapunov, Bohl and Sacker-Sell spectral intervals for
  differential-algebraic equations}.
\newblock {\em Journal of Dynamics and Differential Equations}, 21(1):153--194,
  2009.

\bibitem{Meyer00}
C.D. Meyer.
\newblock {\em {Matrix analysis and applied linear algebra}}.
\newblock Society for Industrial and Applied Mathematics, 2000.

\bibitem{FO09}
F.~Ollivier.
\newblock {Jacobi's bound and normal forms computations. A historical survey}.
\newblock {\em Arxiv preprint arXiv:0911.2674}, 2009.

\bibitem{Pryce01}
J.D. Pryce.
\newblock {A simple structural analysis method for DAEs}.
\newblock {\em BIT Numerical Mathematics}, 41(2):364--394, 2001.

\bibitem{RARH96}
P.J. Rabier and W.C. Rheinboldt.
\newblock {Classical and generalized solutions of time-dependent linear
  differential- algebraic equations* 1}.
\newblock {\em Linear Algebra and its Applications}, 245:259--293, 1996.

\bibitem{RMB00}
G.~Reissig, W.S. Martinson, and P.I. Barton.
\newblock {Differential-algebraic equations of index 1 may have an arbitrarily
  high structural index}.
\newblock {\em SIAM Journal on Scientific Computing}, 21(6):1987--1990, 2000.

\bibitem{rh84}
W.C. Rheinboldt.
\newblock {Differential-algebraic systems as differential equations on
  manifolds}.
\newblock {\em Mathematics of computation}, 43(168):473--482, 1984.

\bibitem{riaza08}
R.~Riaza.
\newblock {\em {Differential-algebraic systems: analytical aspects and circuit
  applications}}.
\newblock World Scientific Pub Co Inc, 2008.

\bibitem{TI08}
M.~Takamatsu and S.~Iwata.
\newblock {Index reduction for differential-algebraic equations by substitution
  method}.
\newblock {\em Linear Algebra and its Applications}, 429(8-9):2268--2277, 2008.

\end{thebibliography}
\end{document}